\documentclass[12pt]{article}
\usepackage{amsfonts}
\usepackage{latexsym,amssymb,amsfonts}
\usepackage{mathrsfs}
\usepackage{graphicx}
\usepackage{psfrag}
\usepackage{amsmath, amsbsy}
\usepackage{amsopn, amstext}
\usepackage{amsmath,multicol,amsopn}
\usepackage{latexsym,amssymb}
\def\sqr#1#2{{\vcenter{\vbox{\hrule height.#2pt
              \hbox{\vrule width.#2pt height#1pt \kern#1pt \vrule width.#2pt}
              \hrule height.#2pt}}}}
\def\signed #1{{\unskip\nobreak\hfil\penalty50
              \hskip2em\hbox{}\nobreak\hfil#1
              \parfillskip=0pt \finalhyphendemerits=0 \par}}
\def\endpf{\signed {$\sqr69$}}

\def\dbR{{\mathop{\rm l\negthinspace R}}}

\def\3n{\negthinspace \negthinspace \negthinspace }
\def\2n{\negthinspace \negthinspace }
\def\1n{\negthinspace }

\def\ds{\displaystyle}

\def\dbR{{\mathop{\rm l\negthinspace R}}}
\def\={\buildrel \triangle \over =}

\def\a{\alpha}
\def\b{\beta}

\def\d{\delta}
\def\e{\varepsilon}

\def\l{\lambda}

 \def\n{\nabla}
\def\si{\sigma}
\def\t{\times}
\def\f{\varphi}

\def\o{\omega}

\def\ns{\noalign{\ss} }

\def\pa{\partial}
\def\G{\Gamma}
\def\D{\Delta}

\def\Si{\Sigma}

\def\O{\Omega}

\def\cF{{\cal F}}

\def\no{\noindent}

\def\ms{\medskip}
\def\bs{\bigskip}
\def\q{\quad}
\def\qq{\qquad}

\def\max{\mathop{\rm max}}

\def\exp{\mathop{\rm exp}}

\def\pa{\partial}

\def\cd{\cdot}
\def\cds{\cdots}

\def\|{\Big |}
\def\({\Big (}
\def\){\Big )}
\def\[{\Big[}
\def\]{\Big]}
\def\be{\begin{equation}}
\def\bel{\begin{equation}\label}
\def\ee{\end{equation}}
\def\bt{\begin{theorem}}
\def\bcd{\begin{condition}}
\def\ecd{\end{condition}}
\def\et{\end{theorem}}
\def\bc{\begin{corollary}}
\def\ec{\end{corollary}}
\def\bde{\begin{definition}}
\def\ede{\end{definition}}
\def\bl{\begin{lemma}}
\def\el{\end{lemma}}
\def\bp{\begin{proposition}}
\def\ep{\end{proposition}}
\def\br{\begin{remark}}
\def\er{\end{remark}}
\def\ba{\begin{array}}
\def\ea{\end{array}}
\def\ed{\end{document}}
\def\ns{\noalign{\ms}}
\def\ds{\displaystyle}

\def\square#1{\vbox{\hrule\hbox{\vrule height#1%
     \kern#1\vrule}\hrule}}
\def\rectangle#1#2{\vbox{\hrule\hbox{\vrule height#1%
     \kern#2\vrule}\hrule}}

\font\tenbb=msbm10 \font\sevenbb=msbm7 \font\fivebb=msbm5

\newfam\bbfam
\scriptscriptfont\bbfam=\fivebb \textfont\bbfam=\tenbb
\scriptfont\bbfam=\sevenbb

\newtheorem{lemma}{Lemma}[section]
\newtheorem{remark}{Remark}[section]

\newtheorem{theorem}{Theorem}[section]
\newtheorem{corollary}{Corollary}[section]

\newtheorem{definition}{Definition}[section]
\newtheorem{proposition}{Proposition}[section]
\newtheorem{condition}{Condition}[section]

\makeatletter
   
   \@addtoreset{equation}{section}
\makeatother

\begin{document}
\title{\bf Carleman Estimate for Stochastic Parabolic Equations and  Inverse Stochastic Parabolic Problems\thanks{This work was  partially supported by the   NSF of China under grants
11101070, and Grant MTM2008-03541 of the MICINN, Spain, Project PI2010-04 of the Basque Government, the ERC Advanced Grant FP7-246775 NUMERIWAVES and the ESF Research Networking Programme OPTPDE.}
}
\author{Qi L\"u \thanks{ Basque Center for Applied Mathematics (BCAM), Mazarredo 14, 48009, Bilbao, Basque Country, Spain; and School of Mathematical Sciences,
University of Electronic Science and Technology
of China, Chengdu 610054, China.{\small\it
E-mail:} {\small\tt luqi59@163.com}.} }

\date{}

\maketitle

\begin{abstract}\no
In this paper, we establish a global Carleman
 estimate for stochastic parabolic equations.
Based on this estimate, we study two inverse
 problems for stochastic parabolic equations.
One  is concerned with a determination problem of the history
of a stochastic heat process through the
observation at the final time $T$, for which we
obtain a conditional stability estimate. The
other is  an inverse source
problem with observation on the lateral
boundary. We derive the uniqueness of the source.
\end{abstract}

\bs

\no{\bf 2010 Mathematics Subject Classification}.  Primary 65N21, 60H15.

\bs

\no{\bf Key Words}. Stochastic parabolic
equations, Carleman estimate, conditional
stability, inverse source problem

\ms

\section{Introduction}

In this paper, we study two different inverse
problems for stochastic parabolic equations  by
 establishing a global Carleman
estimate. We first introduce some notations.

Let $T > 0$, $G \subset \mathbb{R}^{n}$ ($n \in
\mathbb{N}$) be a given bounded domain with a
$C^{2}$ boundary $\G$. Put $$Q \= (0,T) \t G,
\q  \Si \= (0,T) \t \G.$$

Let $(\O, {\cal F}, \{{\cal F}_t\}_{t \geq 0},
P)$ be a complete filtered probability space on
which a  one dimensional standard Brownian motion
$\{ B(t) \}_{t\geq 0}$ is defined. Let $H$ be a
Banach space. Denote by $L^{2}_{\cal F}(0,T;H)$
the Banach space consisting of all $H$-valued $\{
{\cal F}_t \}_{t\geq 0}$-adapted processes
$X(\cdot)$ such that
$\mathbb{E}(|X(\cdot)|^2_{L^2(0,T;H)}) < \infty$,
with the canonical norm; by $L^{\infty}_{\cal
F}(0,T;H)$ the Banach space consisting of all
$H$-valued $\{ {\cal F}_t \}_{t\geq 0}$-adapted
bounded processes; by $L^{2}_{\cal
F}(\O;C([0,T];H))$ the Banach space consisting of
all $H$-valued $\{ {\cal F}_t \}_{t\geq
0}$-adapted processes $X(\cdot)$ satisfying that
$\mathbb{E}(|X(\cdot)|^2_{C(0,T;H)}) < \infty$,
with the canonical norm(similarly, one can define
$L^2_\cF(\O;C^k([0,T];H))$ for any positive $k$).

Throughout this paper, we make the following
assumptions on the coefficients $$b^{ij}:
\;\O\times Q \to  \dbR ,\q
(i,j=1,2,\cdots,n):$$


{\bf (H1)} {\it   $b^{ij}\in
L_{\cF}^2(\O;C^1([0,T];W^{2,\infty}(G)))$ and $b^{ij}=b^{ji}$;}

\ms

 {\bf (H2)} {\it There is a constant $\si>0$ such that
 \bel{h1}
 \sum_{i,j=1}^n b^{ij}(\o,t,x)\xi^{i}\xi^{j}
 \geq \si |\xi|^{2},\q (\o,t,x,\xi)\equiv (\o,t,x,\xi^{1},\cdots,\xi^{n})
 \in \Omega\times Q\t \dbR^{n}.
 \ee}\no
\no

 Let
$$a_1 \in L^{\infty}_{\cF}(0,T;L^{\infty}(G;\dbR^n)), \;\; a_2\in
L^{\infty}_{\cF}(0,T;L^{\infty}(G)) \;\;  \;a_3\in
L^{\infty}_{\cF}(0,T; W^{1,\infty}(G)),$$
 $$
f \in L^2_{\cF}(0,T;L^2(G)) \q\mbox{and }\; g\in L^2_{\cF}(0,T;H^1(G)).
$$
Consider the following stochastic parabolic
equation:
\begin{eqnarray}{\label{system3bu}}
\left\{
\begin{array}{lll}\ds
dy-\sum_{i,j=1}^n (b^{ij}y_{x_i})_{x_j}dt = \big[
(a_1,\nabla y) + a_2 y + f \big] dt + (a_3 y + g
)dB & {\mbox { in }} Q,
 \\
\ns\ds y = 0 & \mbox{ on } \Si,\\
\ns\ds y(0)=y_0, &\mbox{ in } \O,
\end{array}
\right.
\end{eqnarray}
where $y_0\in L^2(\O,\cF_0,P;L^2(G))$ and $y_{x_i} = \frac{\pa y}{\pa x_i}$.

We first recall the definition of the weak and
strong  solution of equation \eqref{system3bu}
and give some  well-posedness results.

\begin{definition}\label{def1}
We call a stochastic process $y\in
L^2_{\cF}(\O;C([0,T];L^2(G)))\cap
L^2_{\cF}(0,T;H_0^1(G))\cap
L^2_{\cF}(\O;C((0,T];H^{1}_0(G)))$ a weak
solution of equation \eqref{system3bu} if for any
$t\in [0,T]$ and any $p\in H_0^1(G)$, it holds
that
\begin{equation}\label{def1eq1}
\begin{array}{ll}
\ds\q\int_G y(t,x)p(x)dx - \int_{G}y_0(x)p(x)dx \\
\ns\ds = \int_0^t \int_G  \Big\{ - \sum_{i,j=1}^n
b^{ij}(s,x) y_{x_i}(s,x)p_{x_j}(x)\\
\ns\ds \qq\qq\;\;
 +  \big[ \big( a_1(s,x), \nabla y(s,x) \big)  +  a_2(s,x) y(s,x) + f(s,x) \big]p(x) \Big\}dxds \\
\ns\ds\q + \int_0^t\int_G  \big[a_3(s,x) y(s,x) +
g(s,x)\big]p(x)dxdB,\qq P\mbox{-a.s.}
\end{array}
\end{equation}
\end{definition}

\begin{definition}\label{def2}
A process $y\in L^2_{\cF}(\O;C([0,T];H^2(G)\cap H_0^1(G)))$ is said to be a strong solution of equation \eqref{system3bu} if for any $t\in [0,T]$, it holds that
\begin{equation}\label{def2eq1}
\begin{array}{ll}
\ds y(t)
 = y_0 + \int_0^t \Big\{-\sum_{i,j=1}^n \big(b^{i,j}(s)y_i(s)\big)_j
  + \big[ (a_1(s),\nabla y(s)) + a_2(s) y(s) + f(s) \big]\Big\}ds \\
\ns\ds\hspace{1.2cm} + \int_0^t\int_G \big[a_3(s)
y(s) + g(s)\big]dB, \qq P\mbox{-a.s.}
\end{array}
\end{equation}
\end{definition}
Obviously, strong solution of equation
\eqref{system3bu} is also its weak solution. We
have the following well-posedness results for
equation \eqref{system3bu}, whose proof  can be
found in \cite[Chapter 6]{Prato}.

\medskip

\begin{lemma}\label{well1}
There exists a unique weak solution of equation \eqref{system3bu}. Furthermore, it holds that
\begin{equation}\label{well1eq1}
|y|_{L^2_\cF(\O;C([0,T];L^2(G)))} + |y|_{L^2_\cF(0,T;H_0^1(G))} \leq Cr_1 |y_0|_{L^2(\O,\cF_0,P;L^2(G))}.
\end{equation}
\end{lemma}
Here and in the sequel, $$r_1\=
|a_1|^2_{L^{\infty}_{\cF}(0,T;L^{\infty}(G;\dbR^n))}
+ |a_2|^2_{L^{\infty}_{\cF}(0,T;L^{\infty}(G))} +
|a_3|^2_{L^{\infty}_{\cF}(0,T;
W^{1,\infty}(G))}+1.$$

\begin{lemma}\label{well2}
Let $y_0\in L^2(\O,\cF_0,P;H^2(G)\cap H_0^1(G))$,
$b_1 \in
L^{\infty}_{\cF}(0,T;W^{1,\infty}(G;\dbR^n))$, $
b_2\in L^{\infty}_{\cF}(0,T;W^{1,\infty}(G))$,
and $b_3\in L^{\infty}_{\cF}(0,T;
W^{2,\infty}(G))$.  Then there exists a unique
strong solution of equation \eqref{system3bu}.
\end{lemma}

Next, we recall the following It\^{o}'s formula,
which plays a key role in the sequel.

\medskip

\begin{lemma}\label{ito}{\rm[{\bf It\^{o}'s formula}]}\q
Let $X(\cd)\in L^2_\cF(0,T;H_0^1(G))$ be a
continuous process with values in $H^{-1}(G)$.
Suppose that there exist $X_0\in
L^2(\O,\cF_0,P;L^2(G))$,  $\Phi(\cd)\in
L^2_\cF(0,T;H^{-1}(G))$ and $\Psi(\cd)\in
L^2_\cF(0,T;L^2(G))$ such that for any $t\in
[0,T]$, it holds that
\begin{equation}\label{ito1}
X(t) = X_0 + \int_0^t \Phi(s)ds + \int_0^t
\Psi(s)dB, P\mbox{-a.s.}
\end{equation}
 in $H^{-1}(G)$.
Then we have that
\begin{equation}\label{ito2}
\begin{array}{ll}\ds
|X(t)|^2_{L^2(G)} = |X(0)|^2_{L^2(G)} +  2\int_0^t \( X(s),\Phi(s) \)_{H_0^{1}(G),H^{-1}(G)}ds\\
\ns\ds\hspace{2.4cm} + 2\int_0^t \( X(s),\Psi(s)
\)_{L^2(G)}dB  + \int_0^t |\Psi(s)|^2_{L^2(G)} ds
\end{array}
\end{equation}
for arbitrary $t\in [0,T]$.
\end{lemma}

\medskip

\begin{remark}
Here we only present a special case for the
It\^{o}'s formula.  It is enough for the proof in
our paper. The general form
  can be found in \cite[Chapter
1]{Pardoux}.
\end{remark}

\begin{remark}
Obviously, both the weak  and strong solution  of
equation \eqref{system3bu} satisfy the
assumptions for Lemma \ref{ito}.  In this paper,
we sometimes use the differential form the the
above It\^{o}'s formula, that is, $d(X^2) = 2XdX
+ (dX)^2$, for the  simplicity of notations.
 \end{remark}

In this paper, we establish a Carleman estimate
for equation \eqref{system3bu}. The so-called
Carleman estimate is a class of weighted energy
estimates which is in connection with (stochastic)
differential operators. As far as we know, the
first example of such kind of estimate appeared
in Carleman's pioneer work for the uniqueness of
the solution of first order elliptic system with
two variables(see \cite{Carleman}).  The idea was
generalized to get the uniqueness of the
solutions for general Cauchy problems in
\cite{Calderon}.  Now it is a useful tool for
studying the uniqueness and unique continuation
property for partial differential equations(see
\cite{Hormander} for example). Such kind of
estimate has been introduced to solving inverse
problems in \cite{Bukhgeim1}, and were
comprehensively studied in
\cite{Isakov,Klibanov1}.  Now it is a helpful
methodology for solving inverse problems (e.g.
\cite{Isakov,Klibanov0,Klibanov1,Yamamoto1,Yamamoto}).
Although the form of Carleman estimate seems to
be very complex, the idea behind them is very
simple. One can understand it by  the following
example.

Let
\begin{equation}\label{ex1}
\left\{
\begin{array}{ll}
\ds \frac{dx}{dt} = a(t)x &\mbox{ in } (0,T],\\
\ns\ds x(0)=x_0.
\end{array}
\right.
\end{equation}
Here $x_0\in \dbR$ and $a(\cd)\in L^\infty(0,T)$.
We prove that there exists a constant $C>0$ such
that for any $x_0\in \dbR$, $|x(T)|\leq C|x_0|$
by Carleman estimate. This result is almost
trivial. And one can prove it without utilizing
Carleman estimate. However, the proof employed
here shows all the ideas of Carleman estimate.

Let $\tilde x(t) = e^{-\varsigma t} x(t)$ with $\varsigma\geq 0$. Then we have
$$
\frac{d \tilde x^2}{dt} = 2a(t)\tilde x^2 - 2\varsigma\tilde x^2 = 2\big[a(t)\ - \varsigma\big] \tilde x^2.
$$
If we choose $\varsigma \geq |a|_{L^\infty(0,T)}$, then we know that $\frac{d \tilde x^2}{dt} \leq 0$, which implies that $\tilde x^2(T) \leq \tilde x^2(0)$. Hence, we get
\begin{equation}\label{exeq1}
e^{-2\varsigma T}x^2(T) \leq x^2(0).
\end{equation}
From this, we obtain  $|x(T)|\leq
e^{T|a|_{L^\infty(0,T)}}|x_0|$ immediately. Thus,
we prove the desired result and we know $C$ can
be chosen to be $e^{T|a|_{L^\infty(0,T)}}$.

Inequality \eqref{exeq1} is a kind of Carleman
estiamte.   The function $e^{-\varsigma t}$ is
called weight function and $\varsigma$ is a
parameter which can be chosen for our purpose. By
means of the choice of $\varsigma$, we control
the lower order term $a(t)x$ and obtain
inequality \eqref{exeq1}. For (stochastic)
partial differential equations, both the choice
of the weight function and the computation are
much more complex. However, they enjoy  the same
idea.

Now we introduce the Carleman estimate to be
established  in this paper. To start with, we
give  some functions. Let $s\in (0,+\infty)$,
$t\in (0,+\infty)$, and $\psi \in
C^{\infty}(\dbR)$ with $|\psi_t|\geq 1$, which is
independent of the $x$-variable. Put
\begin{equation}\label{weight1} \varphi=e^{\l
\psi} \;\mbox{ and } \;\theta = e^{s\varphi}.
\end{equation}

We have the following result.

\medskip

\begin{theorem}\label{carleman est1}
Let $\d\in[0,T)$. Let $\varphi$ and $\theta$ be given in \eqref{weight1}. There exists
a $\l_1>0$ such that for all $\l\geq\l_1$,
there exists an $s_0(\l_1)>0$ so that for all $s\geq
s_0(\l_1)$,  it holds that
\begin{eqnarray}\label{car1}
\begin{array}{ll}\ds
  \q\l \mathbb{E}\int_\d^T\int_G \theta^2|\n y|^2 dxdt + s\l^2 \mathbb{E}\int_\d^T\int_G
\varphi \theta^2 y^2 dxdt \\
\ns\ds \leq C\,\mathbb{E}\Big[ \theta^2(T)|\n y(T)|_{L^2(G)}^2 + \theta^2(\d)|\n
y(\d)|_{L^2(G)}^2 + s\l\varphi(T)\theta^2(T)|y(T)|_{L^2(G)}^2 \\
\ns\ds \q + s\l\varphi(\d)\theta^2(\d)
|y(\d)|_{L^2(G)}^2+\int_\d^T\int_G (1+\f)\theta^2\big( f^2 + g^2 + |\n
g|^2\big)dxdt \Big],
\end{array}
\end{eqnarray}
Here $y$ is arbitrary weak solution of
 equation \eqref{system3bu}.
\end{theorem}

Here and in the sequel, the constant $C$ depends only on $G$,
$(b^{ij})_{n\t n}$, $T$, $\d$ and $\psi$, which may change from line
to line.

Although there are numerous results for the
global Carleman estimate for deterministic
parabolic equations(see
\cite{Fursikov-Imanuvilov1,Yamamoto1} for
example), people know very little about the
stochastic counterpart. In fact, as far as we
know, \cite{barbu1,Tang-Zhang1} are the only two
published papers addressing  the global Carleman
estimate for stochastic parabolic equations.   In
\cite{barbu1,Tang-Zhang1}, some Carleman-type
inequalities  were established, for  deriving the
null controllability of stochastic parabolic
equations.  Note further that the weight function
$\theta$ used in this paper (which plays a key
role in the sequel) is quite different from that
in \cite{barbu1,Tang-Zhang1}. It seems that the
Carleman  estimate in \cite{barbu1,Tang-Zhang1}
cannot be applied to studying the inverse
problems introduced in ths sequel. Indeed, the
weight function $\theta$ in
\cite{barbu1,Tang-Zhang1} is supposed to vanish
at $0$ and $T$, and therefore it does not serve
the purpose of proving Theorem \ref{inv th1} and
Theorem \ref{inv th2}.

As applications of Theorem \ref{carleman est1},
we study  two inverse problems for stochastic
parabolic equations.  There are abundant  works
addressing  the inverse problems for PDEs. And it
is even impossible to list the related papers
owing to the big amount.  However, there exist a
very few works addressing  inverse problems for
stochastic PDEs (see \cite{BCLZ,CT,IK} for
example). Although there are some people
considering the inverse source  problem for
parabolic equations with random noise in the
measurement (see \cite{Johansson1} for example),
to the best of our knowledge, there is no paper
considering the inverse problem for stochastic
parabolic equations.

Now we introduce the inverse problems studied in
this paper. Consider the following stochastic
parabolic equation:
\begin{eqnarray}{\label{system1bu}}
\left\{
\begin{array}{lll}\ds
dy-\sum_{i,j=1}^n(b^{ij}y_{x_i})_{x_j}dt =
[(a_1,\nabla y) + a_2 y ] dt + a_3 y dB  &
{\mbox { in }} Q,
\\
\ns\ds y = 0 & \mbox{ on } \Si,\\
y(0)=y_0 & \mbox{ in } \O.
\end{array}
\right.
\end{eqnarray}
Here $y_0\in L^2(\O,\cF_0,P)$.

\vspace{0.1cm}

The first inverse problem is  concerned with
the following problem:

\vspace{0.2cm}

{\bf Stochastic parabolic equation backward in time:}
 {\it  Let $0 \leq t_0 < T$. Determine $y(\cd,t_0)$, $P$-a.s. from $y(\cd,T)$.}\\

\vspace{0.15cm}

  For deterministic
parabolic equations, such kind of problem has
lots of applications in mathematical physics
(e.g. \cite{Ames}) and is studied extensively
(see \cite{Yamamoto} for a nice survey).
Generally speaking, the problem of   (stochastic)
parabolic equation backward in time is ill-posed.
Small errors in the measuring of the terminal
data may cause huge deviations in final results,
that is, there is no stability in this problem.
Fortunately, if we assume a priori bound for
$y(0)$ (such assumption is reasonable from a
practical viewpoint), then we can regain the
stability in some sense. The concept of
conditional stability is used to describe such
kind of stability.  In general framework, the
conditional stability problem   can be formulated
as follows:

\vspace{0.5cm}

{\it Let $t_0\in [0,T)$, $\a_1\geq 0$, $\a_2\geq
0$  and $M>0$. Put $$U_{M,\a_1} \= \{f \in
L^2(\O,P,\cF_0;H^{\a_1}(G))\,:\,|f|_{
L^2(\O,P,\cF_0;H^{\a_1}(G))} \leq M\}.$$ If
$y_0\in U_{M,\a_1}$,  then can we choose a
function $\b\in C[0,+\infty)$ satisfying the following properties: $$ \left\{  \begin{array}{ll}\ds1. \;\b \geq 0 \mbox{ and } \b \mbox{ is strictly increasing};\\
 \ns\ds 2. \;\lim_{\eta\to 0}\b(\eta)=0;\\
 \ns\ds 3.\; |y(t_0)|_{L^2(\O,\cF_{t_0},P;H^{\a_1}(G))} \leq \b\big(|y(T)|_{L^2(\O,\cF_{T},P;H^{\a_2}(G))}\big).
 \end{array}
 \right.$$    }

\begin{remark}
Here we expect the existence of $\b$ with the
assumptions that $y_0$ belongs to a special set
$U_{M,\a_1}$, which means that $y_0$ enjoys a
priori bound in some sense. Generally speaking,
 $\b$ depends on $M$ and $\a_1$. Once we
choose $M$ and $\a_1$, we add some conditions to
the initial data of equation \eqref{system1bu}.
Hence, that the stability result implied by $\b$
depends on our choice of the initial data. This
is  why we call it ``conditional stability".
\end{remark}

\begin{remark}
The first property for $\b$ means that we only
choose $\b$ in a special class of functions, that
is, the strictly increasing functions. The second
and the third property guarantee  the conditional
stability. Without assuming property 2, we can
always construct $\b$ as $\b(x)=C+x$ with a
constant $C$ which is large enough. However, such
kind of functions do not make any sense for
conditional stability.
\end{remark}

\begin{remark}
Once $\b$ exists, it is not unique. For example,
$\tilde b(x)=\b(x)+x$ is another function
satisfying the three properties.
\end{remark}

In this paper, we obtain the following
interpolation  inequality for the weak solution
of equation \eqref{system1bu}, which implies a
conditional stability result for  equation
\eqref{system1bu} backward in time.

\medskip

\begin{theorem}\label{inv th1}
Let $t_0\in [0,T]$. Then there exist a constant
$\theta\in (0,1)$ and a  constant $C>0$ such that
\begin{equation}\label{conditional sta1}
|y(t_0)|_{L^2(\O,\cF_{t_0},P;L^2(G))} \leq
C\,|y|^{1-\theta}_{L^2_{\cF}(0,T;L^2(G))}
|y(T)|^{\theta}_{L^2(\O,\cF_T,P;H^1(G))},
\end{equation}
for any $y$ solving equation \eqref{system1bu} in
the sense of weak solution.
\end{theorem}

\medskip

As a consequence, we obtain the following result.

\medskip

\begin{theorem}\label{th consta}
Let $y_0\in U_{M,0}$, $\a_2=1$ and
$\b(x)=CM^{1-\theta}x^{\theta}$ with a constant
$C$ independent of $y(0)$. Then we have
$$
|y(t_0)|_{L^2(\O,\cF_{t_0},P;L^{2}(G))} \leq \b\big(|y(T)|_{L^2(\O,\cF_{T},P;H^{1}(G))}\big).
$$
\end{theorem}

\medskip

The proof of Theorem \ref{th consta} follows
Lemma \ref{well1} and Theorem \ref{inv th1}
immediately. We omit it here.

In deterministic setting, a    result which is stronger than Theorem \ref{inv th1}
 was
obtained in \cite{Phung1},  where the authors
study the following equation:
\begin{eqnarray}{\label{dsystem1bu}}
\left\{
\begin{array}{lll}\ds
y_t-\D y = by & {\mbox { in }} Q,
\\
\ns\ds y = 0 & \mbox{ on } \Si.
\end{array}
\right.
\end{eqnarray}
Here $b$ is a suitable function. With the
assumption that $G$ is convex, they get
\begin{equation}\label{dests}
|y(0)|^2_{L^2(G)} \leq
C\exp\left(\frac{|y(0)|_{L^2(G)}}{|y(0)|_{H^{-1}(G)}}\right)|y(T)|^2_{L^2(G_0)}.
\end{equation}
Here $G_0$ is any open subset of $G$. Compared
with Theorem \ref{inv th1}, only
$|y(T)|^2_{L^2(G_0)}$ is involved in the right
hand side of the inequality.    They prove this
result by employing some special frequency
functions, which were first constructed for
proving the doubling property of the solution of
heat equations. However, since the solution of
equation \eqref{system1bu} is non-differentiable
with respect to $t$, it seems that their method
cannot be easily adopted to solve our problem.

\vspace{0.2cm}

As another consequence of Theorem \ref{inv th1}, we
get a backward uniqueness   for  equation \eqref{system1bu}.

\medskip

\begin{corollary}\label{backuni}
Assume that $y $ is a weak solution of equation
(\ref{system1bu}).
If $y(T)=0$ in $G$, $P$-a.s., then $y(t)=0$ in
$G$, $P$-a.s. for all $t\in [0,T]$.
\end{corollary}

The uniqueness problems for the solutions  of
both deterministic and stochastic partial
differential equations have been studied for a
long time. There are a great many positive
results and some negative results.  In case of
time reversible systems, the backward uniqueness
is equivalent to the classical (forward)
uniqueness. If one considers  time irreversible
systems, such as parabolic equations, the
situation is quite different. The backward
uniqueness implies the classical (forward)
uniqueness, however, generally speaking, the
converse conclusion is   untrue.

On
account of the plentiful applications, such as
studying the long time behavior of solutions
and establishing the approximate
controllability from the null controllability,
the backward uniqueness for
parabolic equations draws lots of attention(see
\cite{ESS,Ghidaglia,Kukavica,Lee-Protter,Santo-Prizzi}
and the references cited therein). It is well
understood now. On the contrast,  as far as we
know, \cite{BN} is the only   paper
concerned with backward uniqueness for stochastic parabolic equations in
the literature. In \cite{BN}, the authors
obtained  the backward uniqueness for  semilinear stochastic parabolic equations
with deterministic coefficients. They employed
some deep tools in Stochastic Analysis to
establish the result. However, it seems  that
their method depends on the very fact that the coefficients
are deterministic and one cannot simply mimic
their method to obtain Corollary \ref{backuni},
since the coefficients are random.

The other inverse problem studied in this paper
is about the global uniqueness of an inverse source
problem for stochastic parabolic equations. We
first give a precise formulation of the problem.

\vspace{0.2cm}

Let $x=(x_1,x')\in \dbR^n$ and $x' =
(x_2,\cdots,x_n)\in \dbR^{n-1}$. Consider a special $G$ as  $G=(0,l)\t
G'$, where $G'\subset \dbR^{n-1}$ be a bounded domain
with a $C^2$ boundary. We consider the
following stochastic parabolic equation:
\begin{equation}\label{system2bu}
\left\{
\begin{array}{ll}\ds
dy - \D y = [ (b_1,\n y) + b_2 y + h(t,x')R(t,x) ]dt + b_3 ydB(t) &\mbox{ in } Q,\\
\ns\ds y = 0 & \mbox{ on } \Si,\\
\ns\ds y(0)=0 &\mbox{ in } G.
\end{array}
\right.
\end{equation}
Here
$$b_1 \in L^{\infty}_{\cF}(0,T;W^{1,\infty}(G;\dbR^n)), \;\; b_2\in
L^{\infty}_{\cF}(0,T;W^{1,\infty}(G)), \;\;
\;b_3\in L^{\infty}_{\cF}(0,T;
W^{2,\infty}(G)),$$ and
$$ R\in C^2([0,T]\t \overline G),\q h\in L^2_{\cF}(0,T;H^1(G')).  $$

 The inverse source
problem studied here is   as follows:

\vspace{0.2cm} {\it  Let $R$ be given and
$0<t_0
< T$. Determine the source function $h(t,x')$, $(t,x')
\in (0,t_0)\t G'$,
 by means of the observation of $\ds\frac{\pa y}{\pa \nu}\Big|_{[0,t_0]\t\partial
 G}$.}

\vspace{0.1cm}

Here $\nu = (\nu^1,\cds,\nu^n)\in \dbR^n$ is the outer normal vector of $\G$.

\vspace{0.1cm}

 We have the following uniqueness result
about the above problem.

\begin{theorem}\label{inv th2}
Let
\begin{equation}\label{R}
 |R(t,x)|
\neq 0 \;\mbox{ for all } (t,x)\in [0,t_0]\t
\overline G.
\end{equation}
If $$\frac{\pa y}{\pa \nu} = 0 \mbox{ on }
[0,t_0]\t \pa G, \;\, P\mbox{-a.s.},$$ then
$$
h(t,x') = 0\;\mbox{ for all } (t,x')\in
[0,t_0]\t G', P\mbox{-a.s.}
$$
\end{theorem}

\medskip

\begin{remark}
  One can
follow  the proof of Theorem \ref{inv th2} to
show  that  Theorem \ref{inv th2} also holds when $\D
y$ is substituted by $\ds\sum_{i,j=1}^n
\big(b^{ij}y_i\big)_j$. Here we consider equation
\eqref{system2bu} for the sake of presenting the
key idea in a simple way.

\end{remark}

In practical problems, it is  important to
specify some proper data so that the parameter
to be reconstructed is uniquely identifiable.
In our model, the data utilized is
  the boundary normal derivative of the
solution. This type of inverse problem   is important
in many branches of engineering sciences. For
examples, an accurate estimation of a pollution
source in a river, a determination of magnitude
of groundwater pollution sources.

In the literature, determining a spacewise
dependent   source function  for parabolic
 equations has been considered comprehensively(see \cite{Cannon,Choulli,Isakov,Johansson,Yamamoto}
and the references cited therein). A classical result
for the deterministic setting is as follows.

\vspace{0.2cm}

Consider the following parabolic equation:

\begin{equation}\label{dsystem2bu}
\left\{
\begin{array}{ll}
\ds y_t - \D y = c_1 \n y + c_2 y + R f &\mbox{
in  } Q,\\
\ns\ds y = 0 &\mbox{ on }\Si.
\end{array}
\right.
\end{equation}
Here $c_1$ and $c_2$ are suitable functions on
$Q$. $R\in L^{\infty}(Q)$, $R_t\in
L^{\infty}(Q)$ and $R(t_0,x)\neq 0$ in
$\overline G$ for some $t_0\in (0,T]$. $f\in L^2(G)$ is
independent of $t$. The authors in
\cite{Imanuvilov1} proved the following result:

\vspace{0.2cm}

{\it Assume that $y\in H^{2,1}(Q)$ and
$y_t \in H^{2,1}(Q)$, then there exists a constant
$C> 0$   such that
\begin{equation}\label{dests1}
|f|_{L^2(G)} \leq C\left(|y(t_0)|_{H^2(G)} +
\Big| \frac{\pa y_t}{\pa \nu}
\Big|_{L^2(0,T;L^2(\G_0))}\right),
\end{equation}
where $\G_0$ is any open subset of $\G$.  }

\vspace{0.2cm}

Compared  with Theorem \ref{inv th2},
inequality \eqref{dests1} gives an explicit estimate for the source
term by $|y(t_0)|_{H^2(G)}$ and $\ds\Big| \frac{\pa y_t}{\pa \nu}
\Big|_{L^2(0,T;L^2(\G_0))}$.   A key step in the proof of
equality \eqref{dests1} is to differentiate the solution of
\eqref{dsystem2bu} with respect to $t$. Unfortunately, the solution
of \eqref{system2bu}  does  not enjoy differentiability
with respect to $t$ since the effect of the stochastic noise. However, we can borrow some idea from the proof of inequality \eqref{dests1}. Although it is impossible for us to assume that the solution of equation \eqref{system2bu} is differentiable  with respect to $t$,  we can show that it is differentiable with respect to $x$ with some assumptions(the assumptions in this paper is enough). In this case, we can show the uniqueness of $h$ if $h$ is independent of some $x_i$($i=1,\cds,n$).  Here we suppose that $h$ is independent of $x_1$.

Obviously, both equation \eqref{system1bu} and
equation \eqref{system2bu} are special examples
of equation \eqref{system3bu}.

\medskip

\begin{remark}
As we have pointed out,  the
non-differentiability with respect to the
variable with noise (say, the time variable
considered in this paper) of the solution of a
stochastic PDE usually leads to substantially new
difficulties in the study of inverse problems for
stochastic PDEs. Another trouble for studying the
inverse problem of stochastic PDEs is  that the
usual compactness embedding result does not
remain true for the solution spaces related to
stochastic PDEs. Due to these new difficulties,
 some useful methods for solving inverse
problems for deterministic PDEs (see
\cite{Isakov,Klibanov1} for example) cannot be
used to solve the corresponding inverse
problems in the stochastic setting.
\end{remark}

\medskip

The rest paper is organized as follows. In
Section 2, we prove Theorem \ref{carleman est1}.
Section 3 is addressed  the proof of Theorem
\ref{inv th1}. At last, in Section 4, we give a
proof for Theorem \ref{inv th2}.

\section{Carleman estimate for stochastic parabolic equations}

In this section, we prove Theorem \ref{carleman
est1}.

We first give a weighted  identity, which  plays
an important role in the proof of Theorem
\ref{carleman est1}.

\begin{proposition}\label{identity1}
 Assume that $u$ is an $H^2(\mathbb{R}^n)$-valued
continuous semi-martingale. Put $v=\theta
u$(recall \eqref{weight1} for the definition of
$\theta$). Then we have the following equality:
\begin{equation}\label{bu1}
\begin{array}{ll}\ds
 \q-\theta\Big[ \!\!\sum_{i,j=1}^n (b^{ij}v_{x_i})_{x_j}+s\l\varphi \psi_t v
\Big]\Big[ du -\!\! \sum_{i,j=1}^n (b^{ij}u_{x_i})_{x_j}dt \Big] +
\frac{1}{4}\l\theta v \Big[du - \!\!\sum_{i,j=1}^n
(b^{ij}u_{x_i})_{x_j}dt\Big] \\
\ns\ds= - \sum_{i,j=1}^n \Big(b^{ij}v_{x_i} dv
+ \frac{1}{4}b^{ij}v_{x_i} vdt\Big)_{x_j} +
\frac{1}{2}d\Big( \sum_{i,j=1}^n b^{ij}v_{x_i}
v_{x_j} -
s\l\varphi\psi_t v^2 + \frac{1}{8}\l v^2 \Big) \\
\ns\ds\q - \Big( \frac{1}{2}\sum_{i,j=1}^n
b^{ij}dv_{x_i}
dv_{x_j}+\frac{1}{2}\sum_{i,j=1}^n
b_t^{ij}v_{x_i} v_{x_j} dt - \frac{1}{4}\l
\sum_{i,j=1}^n b^{ij}v_{x_i} v_{x_j} dt \Big) \\
\ns\ds \q + \frac{1}{2}s\l^2\varphi \psi_t^2 v^2 dt + \frac{1}{2}
s\l\varphi\psi_{tt} v^2dt  -
\frac{1}{4}s\l^2\psi_t \varphi v^2 dt + \frac{1}{2}s\l\varphi\psi_t (dv)^2 - \frac{1}{8}\l
(dv)^2 \\
\ns \ds \q   + \Big[ \sum_{i,j=1}^n (b^{ij}v_{x_i})_{x_j}+ s\l\varphi \psi_t v
\Big]^2 dt.
\end{array}
\end{equation}
\end{proposition}

{\it Proof} : The proof is based on some direct
computation by It\^{o}'s stochastic calculus. The first
term in the left hand side of equality
(\ref{bu1}) reads as
\begin{equation}\label{bu2}
\begin{array}{ll}\ds
 \q -\theta\Big[ \sum_{i,j=1}^n
(b^{ij}v_{x_i})_{x_j} + s\l\varphi\psi_t
v \Big] \Big[ du - \sum_{i,j=1}^n (b^{ij}u_{x_i})_{x_j}dt \Big]  \\
 \ns\ds = -\Big[ \sum_{i,j=1}^n
(b^{ij}v_{x_i})_{x_j} + s\l\varphi \psi_t v
\Big]\Big[ dv - \sum_{i,j=1}^n (b^{ij}v_{x_i})_{x_j}dt - s\l\varphi\psi_t v dt \Big] \\
\ns\ds= -\sum_{i,j=1}^n (b^{ij}v_{x_i})_{x_j} dv
-s\l\varphi\psi_t v dv  + \Big[ \sum_{i,j=1}^n
(b^{ij}v_{x_i})_{x_j} +
s\l\varphi \psi_t v \Big]^2 dt\\
\ns\ds = -\sum_{i,j=1}^n (b^{ij}v_{x_i} dv)_{x_j}
+ \frac{1}{2}\sum_{i,j=1}^n d(b^{ij}v_{x_i}
v_{x_j}) -
\frac{1}{2}\sum_{i,j=1}^n b^{ij}dv_{x_i} dv_{x_j}  \\
\ns\ds\q - \frac{1}{2}\sum_{i,j=1}^n b^{ij}_t
v_{x_i} v_{x_j} dt -
\frac{1}{2}d(s\l\varphi\psi_t v^2) \!+\!
\frac{1}{2}s\l^2\psi_t^2
\varphi v^2 dt  \\
\ns\ds \q\,+ \frac{1}{2} s\l\varphi\psi_{tt} v^2
dt +  \frac{1}{2}s\l\psi_t\varphi(dv)^2  + \Big[
\sum_{i,j=1}^n (b^{ij}v_{x_i})_{x_j} + s\l\varphi\psi_t
v  \Big]^2 dt.\nonumber
\end{array}
\end{equation}
The second term in the left hand side of equality (\ref{bu1})
satisfies
\begin{equation}\label{bu3}
\begin{array}{ll}\ds
\q \frac{1}{4}\l\theta v \Big[du - \sum_{i,j=1}^n
(b^{ij}u_{x_i})_{x_j}dt\Big]  \\
\ns\ds  =  \frac{1}{4}\l v\Big[ dv
-\sum_{i,j=1}^n (b^{ij}v_{x_i})_{x_j}dt -
s\l\varphi\psi_t v dt \Big] \\
\ns \ds  =  \frac{1}{4}\l vdv - \frac{1}{4}\l v
\sum_{i,j=1}^n
(b^{ij}v_{x_i})_{x_j}dt - \frac{1}{4}s\l^2\varphi\psi_t v^2 dt \\
\ns \ds  = \frac{1}{8}\l dv^2 - \frac{1}{8}\l
(dv)^2 - \frac{1}{4}\l\sum_{i,j=1}^n
(b^{ij}v_{x_i} v)_{x_j} dt +
\frac{1}{4}\l\sum_{i,j=1}^n b^{ij}v_{x_i} v_{x_j}
dt - \frac{1}{4}s\l^2\psi_t \varphi v^2
dt.\nonumber
\end{array}
\end{equation}
This, together with equality (\ref{bu2}),
implies equality (\ref{bu1}).\endpf

\vspace{0.2cm}

Now we are in a position to prove Theorem
\ref{carleman est1}.

\vspace{0.2cm}

{\it Proof of  Theorem \ref{carleman est1}} :
Applying Proposition \ref{identity1} to
equation (\ref{system3bu}) with $u=y$,
integrating equality (\ref{bu1}) on $[\d,T]\t
G$ for some $\d\in [0,T)$, and taking mathematical
expectation, we get that
\begin{equation}\label{buil1}
\begin{array}{ll}\ds
\q-\mathbb{E}\int_\d^T \int_G \theta\Big[
\sum_{i,j=1}^n (b^{ij}v_{x_i})_{x_j} +
s\l\varphi\psi_t v \Big] \Big[  du  -
\sum_{i,j=1}^n
(b^{ij}u_{x_i})_{x_j}dt  \Big]dx  \\
\ns\ds\q + \frac{1}{4}\l\mathbb{E}\int_\d^T
\int_G \theta v  \Big[ dy -  \sum_{i,j=1}^n
(b^{ij}y_{x_i})_{x_j}dt \Big]dx \\
\ns\ds  = \!-
\mathbb{E}\!\int_\d^T\!\!\int_G\sum_{i,j=1}^n
\Big(\!b^{ij}v_{x_i} dv \!+\! \frac{1}{4}\l
b^{ij}v_{x_i} vdt\!\Big)_{x_j}dx \!+\!
\frac{1}{2}\mathbb{E}\!\!\int_\d^T\!\!\int_G \!d\Big(\!
\sum_{i,j=1}^n b^{ij}v_{x_i} v_{x_j} \!-\!
s\l\varphi\psi_t v^2 \!+\! \frac{1}{8}\l v^2 \!\Big)dx\\
\ns\ds \q  -
\mathbb{E}\int_\d^T\int_G\Big(\frac{1}{4}\l
\sum_{i,j=1}^n b^{ij}v_{x_i} v_{x_j} dt +
\frac{1}{2}\sum_{i,j=1}^n b_t^{ij}v_{x_i}
v_{x_j} dt - \frac{1}{2}\sum_{i,j=1}^n
b^{ij}dv_{x_i}
dv_{x_j}   \Big)dx \\
\ns\ds  \q + \mathbb{E}\int_\d^T\int_G\Big[\frac{1}{2}s\l^2\varphi\psi_t^2 v^2 dt +  \frac{1}{2} s\l\varphi\psi_{tt} v^2dt - \frac{1}{4}s\l^2\varphi\psi_t v^2 dt +
\frac{1}{2}s\l\varphi\psi_t(dv)^2 - \frac{1}{8}\l (dv)^2\Big]dx\\
\ns\ds\q + \mathbb{E}\int_\d^T\int_G\Big[
\sum_{i,j=1}^n (b^{ij}v_{x_i})_{x_j} +
s\l\varphi \psi_t v \Big]^2 dxdt.
\end{array}
\end{equation}

Now we estimate the terms in the right  hand
side of equality \eqref{buil1}  one by one.

\vspace{0.2cm}

 For the first one, since $y|_\Si = 0$, we have that $v|_\Si=0$. Therefore, it holds that
\begin{equation}\label{buil2}
\begin{array}{ll}\ds
\q- \mathbb{E}\int_\d^T\int_G\sum_{i,j}^n
\Big(b^{ij}v_{x_i} dv + \frac{1}{4}\l
b^{ij}v_{x_i} vdt\Big)_{x_j}dx\\
\ns\ds = -
\mathbb{E}\int_\d^T\int_{\G}\sum_{i,j=1}^n
b^{ij}\Big( v_{x_i} dv + \frac{1}{4}\l v_{x_i}
vdt\Big)\nu^j d\,\G = 0.
\end{array}
\end{equation}

For the second one,  we have
\begin{eqnarray}\label{buil3}
\begin{array}{ll}
\ds\q\frac{1}{2}\mathbb{E}\int_\d^T\int_G
d\Big( \sum_{i,j=1}^n b^{ij}v_{x_i} v_{x_j} -
s\l\varphi \psi_t v^2 + \frac{1}{8}\l v^2
\Big)dx \\
\ns\ds\geq -C\mathbb{E}\Big(|\n v(T)|_{L^2(G)}^2 + |\n
v(\d)|_{L^2(G)}^2 + s\l\varphi(T)|v(T)|_{L^2(G)}^2 + s\l\varphi(\d)
|v(\d)|_{L^2(G)}^2\Big).
\end{array}
\end{eqnarray}

Since
$$
\begin{array}{ll}\ds
\q\mathbb{E}\int_\d^T\int_G
\frac{1}{2}\sum_{i,j=1}^n b^{ij}dv_{x_i} dv_{x_j}
dx \\
\ns\ds = \frac{1}{2}\mathbb{E}\int_\d^T\int_G
\sum_{i,j=1}^n b^{ij}\theta^2\big(a_3y+g
\big)_{x_i} \big(a_3y+g \big)_{x_j} dxdt,
\end{array}
$$
 the third one reads as
\begin{eqnarray}\label{buil4}
\begin{array}{ll}
\ds\q\mathbb{E}\int_\d^T\int_G\Big(\frac{1}{4}\l
\sum_{i,j=1}^n b^{ij}v_{x_i} v_{x_j} dt +
\frac{1}{2}\sum_{i,j=1}^n b_t^{ij}v_{x_i}
v_{x_j} dt - \frac{1}{2}\sum_{i,j=1}^n
b^{ij}dv_{x_i} dv_{x_j}
\Big)dx \\
\ns\ds \geq \mathbb{E}\int_\d^T\int_G\Big[\frac{1}{4}\l \si |\n v|^2  -C|\n
v|^2  - C \big(a_3^2 |\n v|^2 + |\n a_3|^2 v^2 + \theta^2|\n g|^2 + \theta^2|g|^2\big)  \Big]dxdt \\
\ns\ds \geq \frac{1}{4}\l
\mathbb{E}\int_\d^T\int_G\si |\n v|^2 dxdt
-C(|a_3|^2_{L^{\infty}_\cF(0,T;W^{1,\infty}(G))}+1)\mathbb{E}\int_\d^T\int_G
(|\n v|^2 + v^2) dxdt\\
\ns\ds \q - C\mathbb{E}\int_\d^T\int_G
\theta^2(|\n g|^2 + g^2) dxdt.
\end{array}
\end{eqnarray}

For the forth one, recalling that $|\psi_t|\geq
1$ and utilizing that
$$
\begin{array}{ll}\ds
\q\mathbb{E}\int_\d^T\int_G\Big[
\frac{1}{2}s\l\varphi\psi_t(dv)^2 - \frac{1}{8}\l
(dv)^2\Big]dx \\
\ns\ds = \mathbb{E}\int_\d^T\int_G\theta^2\Big[
\frac{1}{2}s\l\varphi\psi_t(a_3 y + g)^2 -
\frac{1}{8}\l (a_3 y + g)^2\Big]dxdt,
\end{array}
$$
we see
\begin{eqnarray}\label{buil5}
\begin{array}{ll}
\ds\mathbb{E}\int_\d^T\int_G\Big[\frac{1}{2}s\l^2\varphi\psi_t^2 v^2 dt +  \frac{1}{2} s\l\varphi\psi_{tt} v^2dt - \frac{1}{4}s\l^2\varphi\psi_t v^2 dt +
\frac{1}{2}s\l\varphi\psi_t(dv)^2 - \frac{1}{8}\l (dv)^2\Big]dx \\
\ns\ds\geq \frac{1}{4}s\l^2\mathbb{E}\int_\d^T\int_G \varphi v^2 dxdt + sO(\l)\mathbb{E}\int_\d^T\int_G \varphi v^2 dxdt - Cs\l\mathbb{E}\int_\d^T\int_G (1+\varphi) \theta^2 g^2 dxdt.
\end{array}
\end{eqnarray}
Thus, we know that there exists a $\l_0>0$ such that for all $\l\geq\l_0$, it holds that
\begin{eqnarray}\label{buil5.1}
\begin{array}{ll}
\ds\mathbb{E}\int_\d^T\int_G\Big[\frac{1}{2}s\l^2\varphi\psi_t^2 v^2 dt +  \frac{1}{2} s\l\varphi\psi_{tt} v^2dt - \frac{1}{4}s\l^2\varphi\psi_t v^2 dt +
\frac{1}{2}s\l\varphi\psi_t(dv)^2 - \frac{1}{8}\l (dv)^2\Big]dx \\
\ns\ds\geq \frac{1}{8}s\l^2\mathbb{E}\int_\d^T\int_G \varphi v^2 dxdt- Cs\l\mathbb{E}\int_\d^T\int_G (1+\varphi) \theta^2 g^2 dxdt.
\end{array}
\end{eqnarray}

Now, we estimate the terms in the left hand side
one by one. By equation (\ref{system3bu}) and
noting that
$$
-\mathbb{E}\int_\d^T\int_G\theta\Big[
\sum_{i,j=1}^n (b^{ij}v_{x_i})_{x_j}+s\l\varphi \psi_t v
\Big] (a_3 y + g )dB dx = 0,
$$
 we know that
\begin{eqnarray}\label{buir1}
\begin{array}{ll}\ds
\q-\mathbb{E}\int_\d^T\int_G\theta\Big[
\sum_{i,j=1}^n (b^{ij}v_{x_i})_{x_j}+s\l\varphi  \psi_t v
\Big]\Big[ dy - \sum_{i,j=1}^n
(b^{ij}y_{x_i})_{x_j}dt \Big]dx \\
\ns\ds = -\mathbb{E}\int_\d^T\int_G\theta\Big[
\sum_{i,j=1}^n (b^{ij}v_{x_i})_{x_j}+s\l\varphi \psi_t  v
\Big]\Big[
(a_1,\nabla y) + a_2 y + f \Big]dtdx\\
\ns\ds \leq \mathbb{E}\int_\d^T\int_G\Big[
\sum_{i,j=1}^n (b^{ij}v_{x_i})_{x_j}+s\l\varphi  \psi_t v
\Big]^2dxdt + \mathbb{E}\int_\d^T\int_G\theta^2\Big[ (a_1,\nabla y) + a_2 y + f \Big]^2dtdx\\
\ns\ds\leq \mathbb{E}\int_\d^T\int_G \Big[
\sum_{i,j=1}^n (b^{ij}v_{x_i})_{x_j}+s\l\varphi \psi_t  v
\Big]^2 dxdt +
3\mathbb{E}\int_\d^T\int_G\theta^2\big(|a_1|^2 |\n u|^2 + a_2^2 u^2 + f^2\big)dxdt  \\
\ns\ds\leq \mathbb{E}\int_\d^T\int_G \Big[
\sum_{i,j=1}^n (b^{ij}v_{x_i})_{x_j}+s\l\varphi \psi_t 
v \Big]^2 dxdt +
3|a_1|^2_{L^\infty_\cF(0,T;L^\infty(G;\dbR^n))}\mathbb{E}\int_\d^T\int_G
|\n v|^2dxdt  \\
\ns\ds \q +
3|a_2|^2_{L^\infty_\cF(0,T;L^\infty(G))}\mathbb{E}\int_\d^T\int_G
 v^2 dxdt + 3\mathbb{E}\int_\d^T\int_G
\theta^2 f^2 dxdt,
 \end{array}
\end{eqnarray}
and that
\begin{eqnarray}\label{buir2}
\begin{array}{ll}
\ds\q \frac{1}{4}\l \mathbb{E}\int_\d^T\int_G \theta v \Big[ du -
\sum_{i,j=1}^n (b^{ij}u_{x_i})_{x_j}dt \Big]dx \\
\ns\ds = \frac{1}{4}\l \mathbb{E}\int_\d^T\int_G
\theta v \Big[
(a_1,\nabla y) + a_2 y + f \Big]dtdx\\
\ns\ds \leq \frac{1}{64}\l
\mathbb{E}\int_\d^T\int_G  v^2 dxdt +
\mathbb{E}\int_\d^T\int_G \theta^2\Big[
(a_1,\nabla y) + a_2 y + f \Big]^2 dtdx\\
 \ns\ds \leq
\frac{1}{64}\l^2 \mathbb{E}\int_\d^T\int_G  v^2
dxdt + 3\mathbb{E}\int_\d^T\int_G
\theta^2\big(|a_1|^2 |\n u|^2 + a_2^2
u^2+f^2\big)dxdt \\
\ns\ds \leq \frac{1}{64}\l^2
\mathbb{E}\int_\d^T\int_G v^2 dxdt +
3|a_1|^2_{L^\infty_\cF(0,T;L^\infty(G;\dbR^n))}\mathbb{E}\int_\d^T\int_G
|\n v|^2dxdt \\
\ns\ds \q +
3|a_2|^2_{L^\infty_\cF(0,T;L^\infty(G))}\mathbb{E}\int_\d^T\int_G
 v^2 dxdt+ 3\mathbb{E}\int_\d^T\int_G
\theta^2 f^2 dxdt.
 \end{array}
\end{eqnarray}

From (\ref{buil1})--(\ref{buir2}), we find
\begin{equation}\label{bui1}
\begin{array}{ll}\ds
\q\frac{1}{4}\l\mathbb{E}\int_\d^T\!\!\int_G
\!\!|\n v|^2 dxdt \!-\! C
\big(|a_1|^2_{L^\infty_\cF(0,T;L^\infty(G;\dbR^n))}
\!\!+\!\!
|a_3|^2_{L^{\infty}_\cF(0,T;W^{1,\infty}(G))}\!+\!1\big)\mathbb{E}\int_\d^T\!\!\int_G
\!\!|\n
v|^2dxdt \\
\ns\ds  \q
 - C
(|a_2|^2_{L^\infty_\cF(0,T;L^\infty(G))} +
|a_3|^2_{L^{\infty}_\cF(0,T;W^{1,\infty}(G))}+1)\mathbb{E}\int_\d^T\int_G
|\n
v|^2dxdt\\
\ns\ds \q+ \Big(\frac{1}{8}s\l^2-
\frac{1}{64}\l^2\Big)\mathbb{E}\int_\d^T\int_G \varphi v^2 dxdt\\
\ns\ds  \leq C\mathbb{E}\Big[|\n v(T)|_{L^2(G)}^2 + |\n
v(\d)|_{L^2(G)}^2 + s\l\varphi(T)|v(T)|_{L^2(G)}^2 + s\l\varphi(\d)
|v(\d)|_{L^2(G)}^2\\
\ns\ds \qq\q +s\l\int_\d^T\int_G (1+\f)\theta^2\big( f^2 + g^2 + |\n g|^2\big)dxdt \Big].
\end{array}
\end{equation}
Recalling that $$r_1 =   |a_1|^2_{L^\infty_\cF(0,T;L^\infty(G;\dbR^n))} +
|a_2|^2_{L^\infty_\cF(0,T;L^\infty(G))} +
|a_3|^2_{L^{\infty}_\cF(0,T;W^{1,\infty}(G))}+1,$$
from inequality \eqref{bui1}, we know that there exists a
$\l_1\geq \max\big\{Cr_1,\;\l_0\big\}$ such
that for all $\l\geq\l_1$, there exists a
$s_0(\l_1)>0$ so that for all $s\geq s_0(\l_1)$, it holds
that
\begin{equation}\label{bui2}
\begin{array}{ll}\ds
 \q\l \mathbb{E}\int_\d^T\int_G |\n v|^2 dxdt + s\l^2 \mathbb{E}\int_\d^T\int_G
\varphi v^2 dxdt  \\
\ns\ds \leq C\mathbb{E}\Big[|\n v(T)|_{L^2(G)}^2 + |\n
v(\d)|_{L^2(G)}^2 + s\l\varphi(T)|v(T)|_{L^2(G)}^2 + s\l\varphi(\d)
|v(\d)|_{L^2(G)}^2\\
\ns\ds \q +s\l\int_\d^T\int_G (1+\varphi)\theta^2\big( f^2 + g^2 + |\n g|^2\big)dxdt
\Big],
\end{array}
\end{equation}
 which implies inequality \eqref{car1} immediately.
 \endpf

\section{Proof for Theorem \ref{inv th1}}

This section is devoted to the proof of Theorem
\ref{inv th1}. We borrow some ideas from
\cite{Yamamoto}.

\vspace{0.1cm}

{\it Proof of Theorem \ref{inv th1}}\,: Choose $t_1$ and $t_2$ such
that $0< t_1 < t_2 < t_0$. Set $\a_k = e^{\l t_k}$
($k=0,1,2$). Let $\rho\in C^\infty(\mathbb{R})$
such that $0\leq\rho\leq 1$ and that
\begin{equation}\label{chi}
\rho = \left\{
\begin{array}{ll}
\ds 1, & t\geq t_2,\\
\ns\ds 0, & t \leq t_1.
\end{array}
\right.
\end{equation}
Let $z=\rho y$,  by means of $y$ solves equation \eqref{system1bu}, we know that $z$ solves
\begin{equation}\label{inv equ1}
\left\{
\begin{array}{ll}
\ds dz - \sum_{i,j=1}^n (b^{ij}z_{x_i})_{x_j} dt=
\big[(a_1,\n z) + a_2 z + \rho_t(t)y\big]dt +
a_3 z dB(t) &\mbox{ in } Q,\\
\ns\ds z=0 &\mbox{ on } \Si,\\
\ns\ds z(0) = 0 &\mbox{  in  }G.
\end{array}
\right.
\end{equation}
Applying Theorem \ref{carleman est1} with
$\psi=t$ and $\d=0$ to equation \eqref{inv
equ1}, for  $\l\geq\l_1$ and $s\geq s_0(\l_1)$, we
have
\begin{equation}\label{inv eq1}
\begin{array}{lll}\ds
\q\l\mathbb{E}\int_Q \theta^2|\n z|^2dxdt +
s\l^2\mathbb{E}\int_Q \theta^2\varphi
|z|^2dxdt\\
\ns\ds \leq C\mathbb{E}\Big[ \theta^2(T)\big|\n
z(T)\big|^2_{L^2(G)} +
s\l\varphi(T)\theta^2(T)\big|z(T)\big|^2_{L^2(G)} +
 \int_Q \theta^2|\rho_t(t)y|^2 dxdt
\Big].
\end{array}
\end{equation}
From the choice of $\rho$, we see that
\begin{equation}\label{inv eq2}
\mathbb{E}\int_Q \theta^2 |\rho_t(t)|^2 y^2dxdt
\leq C \int_{t_1}^{t_2}\int_G \theta^2 y^2 dxdt
\leq C
\theta^2(t_1)|y|^2_{L^2_{\cF}(0,T;L^2(G))}.
\end{equation}
This, together with inequality \eqref{inv eq1},
implies that
\begin{equation}\label{inv eq3}
\begin{array}{lll}\ds
\q\l\theta^2(t_0)\mathbb{E}\int_{t_0}^T\int_G
 |\n y|^2dxdt +
s\l^2\theta^2(t_0)\mathbb{E}\int_{t_0}^T\int_G
\varphi
|y|^2dxdt\\
\ns\ds \leq \l\mathbb{E}\int_Q \theta^2|\n z|^2dxdt +
s\l^2\mathbb{E}\int_Q \theta^2\varphi
|z|^2dxdt\\
\ns\ds \leq
C\theta^2(t_1)|y|^2_{L^2_{\cF}(0,T;L^2(G))} +
C\mathbb{E}\Big( \theta^2(T)\big|\n
y(T)\big|^2_{L^2(G)} +
s\l\varphi(T)\theta^2(T)\big|y(T)\big|^2_{L^2(G)}
\Big).
\end{array}
\end{equation}
Here we utilize the fact that
$\theta(t)\leq\theta(s)$  for $t\leq s$.

From inequality \eqref{inv eq3}, we see
\begin{equation}\label{inv eq3.1}
\begin{array}{lll}\ds
\q\l\mathbb{E}\int_{t_0}^T\int_G
 |\n y|^2dxdt +
s\l^2\mathbb{E}\int_{t_0}^T\int_G
\varphi
|y|^2dxdt\\
\ns\ds \leq
C\theta^2(t_1)\theta^{-2}(t_0)|y|^2_{L^2_{\cF}(0,T;L^2(G))} +
C\mathbb{E}\Big( \theta^2(T)\big|\n
y(T)\big|^2_{L^2(G)} +
s\l\varphi(T)\theta^2(T)\big|y(T)\big|^2_{L^2(G)}
\Big).
\end{array}
\end{equation}

By means of $ d(y^2) = 2ydy + (dy)^2$, we obtain
that
\begin{equation}\label{inv eq4}
\begin{array}{ll}\ds
\q\mathbb{E}\int_G |y(t_0)|^2dx \\\ns\ds=
\mathbb{E}\int_G |y(T)|^2dx -
\mathbb{E}\int_{t_0}^T\int_G \big[2ydy +
(dy)^2\big]dx\\
\ns\ds =\mathbb{E}\int_G |y(T)|^2dx -
\mathbb{E}\int_{t_0}^T\int_G
\Big\{2y\Big[\sum_{i,j=1}^n
(b^{ij}y_{x_i})_{x_j} +
(a_1,\n y) + a_2 y \Big]+ (a_3 y)^2 \Big\}dxdt\\
\ns\ds \leq \mathbb{E}\int_G |y(T)|^2dx +
C\mathbb{E}\int_{t_0}^T\int_G |\n y|^2 dxdt\\
\ns\ds\q +
\big(|a_1|^2_{L^\infty_\cF(0,T;L^\infty(G;\dbR^n))}
+
|a_2|_{L^\infty_\cF(0,T;L^\infty(G))}+|a_2|^2_{L^\infty_\cF(0,T;L^\infty(G))}\big)\mathbb{E}\int_{t_0}^T\int_G
y^2 dxdt\\
\ns\ds \leq \mathbb{E}\int_G |y(T)|^2dx +
C\mathbb{E}\int_{t_0}^T\int_G |\n y|^2 dxdt +
Cr_1 \mathbb{E}\int_{t_0}^T\int_G y^2 dxdt,
\end{array}
\end{equation}
Recalling $\varphi\geq 1$, from inequality \eqref{inv eq4},  we know that there exists a $\l_2>0$
such that for all $\l\geq \l_2$, it holds that
\begin{equation}\label{inv eq4.1}
\begin{array}{ll}
\ds \q\mathbb{E}\int_G |y(t_0)|^2dx \\\ns\ds
\leq \mathbb{E}\int_G |y(T)|^2dx +
C\Big(\l \mathbb{E}\int_{t_0}^T\int_G
|\n y|^2dxdt +
s\l^2 \mathbb{E}\int_{t_0}^T\int_G
\varphi y^2dxdt \Big).
\end{array}
\end{equation}
Combing   inequality \eqref{inv eq3.1} and
inequality \eqref{inv eq4.1}, for any   $\l\geq
\max\{\l_1,\l_2\}$ and $s\geq s_0(\l_1)$, we have
\begin{equation}\label{inv eq5}
\begin{array}{ll}\ds
\q\mathbb{E}\int_G |y(t_0)|^2dx \\
\ns\ds \leq C
\theta^2(t_1)\theta^{-2}(t_0)\big|y\big|^2_{L^2_{\cF}(0,T;L^2(G))}
+ C\mathbb{E}\left( \theta^2(T)\big|\n
y(T)\big|^2_{L^2(G)} +
s\l\varphi(T)\theta^2(T)\big|y(T)\big|^2_{L^2(G)}
\right).
\end{array}
\end{equation}
 Now we fix $\l_3=\max\{\l_1,\l_2\}$, from
inequality \eqref{inv eq5}, we get
\begin{equation}\label{inv eq6}
\mathbb{E}\int_G |y(t_0)|^2dx \leq C
\theta^2(t_1)\theta^{-2}(t_0)\big|y\big|^2_{L^2_{\cF}(0,T;L^2(G))}
+ C\,\theta^2(T)\mathbb{E}\big| y(T)\big|^2_{H^1(G)}.
\end{equation}
Replacing $C$ by $Ce^{s_0e^{\l_3 T}}$, from
inequality \eqref{inv eq6}, for any $s>0$, it
holds that
\begin{equation}\label{inv eq7}
\mathbb{E}\int_G |y(t_0)|^2dx \leq
Ce^{-2s(e^{\l_3 t_1}-e^{\l_3
t_0})}\big|y\big|^2_{L^2_{\cF}(0,T;L^2(G))} +
Ce^{Cs}\mathbb{E}\big| y(T)\big|^2_{H^1(G)}.
\end{equation}
 Choosing $s \geq 0$ which minimize
the right-hand side of inequality \eqref{inv
eq7}, we obtain that
\begin{equation}\label{inv eq8}
\mathbb{E}\big|y(t_0)\big|^2_{L^2(G)} \leq C
|y|^{1-\theta}_{L^2_{\cF}(0,T;L^2(G))}
\mathbb{E}\big|y(T)\big|^{\theta}_{H^1(G)},
\end{equation}
with
$$
\theta = \frac{2(e^{\l_3 t_0}-e^{\l_3
t_1})}{C+2(e^{\l_3 t_0}-e^{\l_3
t_1})}.
$$
\endpf


\section{Proof of Theorem \ref{inv th2}}

Thie section is devoted to   proving Theorem
\ref{inv th2}. We borrow some ideas in
\cite{Yamamoto} again .

\vspace{0.1cm}

{\it Proof of Theorem \ref{inv th2}}\,:  From the assumptions on $b_1$, $b_2$, $b_3$, $R$ and $h$, and by Lemma \ref{well2}, we know equation \eqref{system2bu} admits a unique strong solution. For arbitrary small $\e>0$, we
choose $t_1$ and $t_2$ such that
$$
0 < t_0 -\e < t_1 < t_2 < t_0.
$$

Let $\chi\in C^{\infty}(\dbR)$ be a cut-off function such
 that $0\leq \chi \leq 1$ and that
\begin{equation}\label{chi1}
\chi = \left\{
\begin{array}{ll}
\ds 1, & t\leq t_1,\\
\ns\ds 0, & t\geq t_2.
\end{array}
\right.
\end{equation}

Put $y = Rz$ (recall \eqref{R} for $R$) in
$[0,t_2]\t G$. Since $y$ is a strong  solution of
equation \eqref{system2bu}, we know that $z$
solves
\begin{equation}\label{system2bu1}
\left\{
\begin{array}{lll}\ds
dz  - \D zdt = \Big[  (b_1, \n z) +  \Big(  \frac{2\n R}{R},  \n z \Big)  +  \Big( b_2  + \frac{\D R}{R}  - \frac{2(\n R,\n R)}{R^2} \\ \ns\ds\hspace{2.5cm}-  \frac{R_t}{R} + \Big( \frac{\n R}{R}, b_1 \Big) \Big)  z\Big]dt  +   h dt
 + b_3 z dB(t) &\mbox{ in } [0,t_0]\t G,\\
\ns\ds z=\frac{\pa z}{\pa \nu} = 0 &\mbox{ on } [0,t_0]\t\G,\\
\ns\ds z(0) = 0 &\mbox{ in } G.
\end{array}
\right.
\end{equation}
Setting $u=z_{x_1}$, noting $z$ is the strong
solution of  equation \eqref{system2bu1}  and
$z_{x_1} = \frac{\pa z}{\pa\nu}=0$ on
$\big(\{0\}\times G'\big)\cup \big(\{l\}\times
G'\big)$, we know that $u$ is the weak solution
of the following equation:

\begin{equation}\label{system2bu2}
\left\{
\begin{array}{lll}\ds
du - \D udt =  \Big[  ((b_1)_{x_1}, \n z) + (b_1,\n u) +  \Big( \Big(\frac{2\n R}{R}\Big)_{x_1},  \n z \Big)  +\Big( \frac{2\n R}{R},  \n u \Big)\\
\ns\ds \hspace{2.5cm} +  \Big(  b_2  +  \frac{\D R}{R}  - \frac{2(\n R,\n R)}{R^2} -  \frac{R_t}{R} +\Big( \frac{\n R}{R},  b_1  \Big)  \Big)_{x_1}  z  \\
\ns\ds\hspace{2.5cm} +\Big( b_2 +  \frac{\D R}{R} - \frac{2(\n R,\n R)}{R^2} - \frac{R_t}{R} + \Big(  \frac{\n R}{R},  b_1  \Big) \Big)  u \Big]dt \\
\ns\ds\hspace{2.5cm} + (b_3)_{x_1} zdt + b_3 u dB(t) &\mbox{ in } [0,t_0]\!\t\! G,\\
\ns\ds u=0 &\mbox{ on } [0,t_0]\!\t\!\G\\
\ns\ds u(0)=0 &\mbox{ in }  G.
\end{array}
\right.
\end{equation}
Set $w = \chi u$. Then we know that $w$ is a weak
solution of the following equation:
\begin{equation}\label{system2bu3}
\left\{
\begin{array}{lll}\ds
 dw \! - \! \D wdt = \! \Big[  ((b_1)_{x_1}, \chi\n z) + (b_1,\n w) +  \Big( \Big(\frac{2\n R}{R}\Big)_{x_1},  \chi\n z \Big)  +\Big( \frac{2\n R}{R},  \n w \Big)\\
\ns\ds \hspace{2.2cm}  + \Big(  b_2  +  \frac{\D R}{R}  - \frac{2(\n R,\n R)}{R^2}-  \frac{R_t}{R}  + \Big( \frac{\n R}{R},   b_1 \Big)  \Big)_{x_1}  \chi z \\
\ns\ds \hspace{2.2cm}+  \Big( b_2 + \frac{\D R}{R}  - \frac{2(\n R,\n R)}{R^2}- \frac{R_t}{R}  +  \Big( \frac{\n R}{R},   b_1 \!\Big) \Big)  w \Big]dt\\
\ns\ds\hspace{2.2cm}   + (b_3)_{x_1}\chi zdB(t) + b_3 w dB(t) - \chi' udt &\mbox{ in } [0,t_0]\!\t\! G,\\
\ns\ds w=0 &\mbox{ on } [0,t_0]\!\t\!\G,\\
\ns\ds w=0 &\mbox{ in } G.
\end{array}
\right.
\end{equation}
By means of  $u = z_{x_i}$ and $z(t,0,x')=y(t,0,x')=0$ for $(t,x')\in (0,t_0)\t G'$, we see
\begin{equation}\label{inv th2 eq1}
\chi z = \chi\int_0^{x_1}u(t,\eta,x')d\eta = \int_0^{x_1} w(t,\eta,x')d\eta.
\end{equation}
This, together with equation \eqref{system2bu3},
implies  that $w$ is the weak solution of the
following equation:
\begin{equation}\label{system2bu4}
\left\{
\begin{array}{lll}\ds
\!\!\!\!dw \!-\! \D wdt \!= \! \Big[  (b_1,\n w) \!+\! \Big( \frac{2\n R}{R},  \n w \Big)  + \Big((b_1)_{x_1}, \n\!\!\! \int_0^{x_1}\!w(t,\eta,x')d\eta\Big) \\
\ns\ds \hspace{2cm}  + \Big( \Big(\frac{2\n R}{R}\Big)_{x_1}, \n \int_0^{x_1}  w(t,\eta,x')d\eta \Big) \\
\ns\ds \hspace{2.1cm} + \Big(  b_2 +  \frac{\D R}{R}  - \frac{2(\n R,\n R)}{R^2}- \frac{R_t}{R}  + \Big(  \frac{\n R}{R},  b_1  \Big) \Big)  w\\
\ns\ds \hspace{2.1cm}  + \Big(  b_2  +  \frac{\D R}{R}  -  \frac{R_t}{R}  +  \Big( \frac{\n R}{R},   b_1  \Big) \!\Big)_{x_1} \int_0^{x_1}w(t,\eta,x')d\eta   \Big]dt\\
\ns\ds\hspace{2.1cm}   + (b_3)_{x_1}\chi \int_0^{x_1}u(t,\eta,x')d\eta dB(t) + b_3 w dB(t) - \chi' udt \!\!\!\!&\mbox{ in } [0,t_0]\!\t\! G,\\
\ns\ds w=0 &\mbox{ on } [0,t_0]\!\t\!\G,\\
\ns\ds w=0 &\mbox{ in } G.
\end{array}
\right.
\end{equation}

Applying Theorem \ref{carleman est1} to equation
\eqref{system2bu4}  with $\psi(t)=-t$, noting
that $w(0) = 0$, and that $w(t_0)=\chi(t_0)u(t_0)
= 0$, we get
\begin{equation}\label{inv th2 ine1}
\begin{array}{lll}\ds
\q\mathbb{E}\int_0^{t_0}\int_G \theta^2 \big(\l |\n w|^2 + s\l^2 w^2  \big)dxdt \\
\ns\ds \leq
C\mathbb{E}\int_0^{t_0}\int_{G}\theta^{2}
|\chi' u|^2 dxdt \\ \ns\ds\q  + Cr_1\mathbb{E}
\int_0^{t_0}\int_G \theta^2 \Big(  \Big|
\int_0^{x_1}w(t,\eta,x')d\eta  \Big|^2 + \Big|
\int_0^{x_1} |\n w(t,\eta,x')d\eta|  \Big|^2
\Big)dxdt.
\end{array}
\end{equation}
Since
$$
\Big| \int_0^{x_1} w(t,\eta,x')d\eta  \Big|^2 \leq l \int_0^l | w(t,\eta,x') |^2d\eta,
$$
we know
\begin{equation}\label{inv th2 ine2}
\begin{array}{ll}\ds
\int_0^{t_0} \int_G \theta^2 \Big| \int_0^{x_1} w(t,\eta,x')d\eta  \Big|^2 dxdt \leq l \int_0^l dx_1 \int_0^{t_0} \int_{G'} \int_0^l \theta^2 |w(t,\eta,x')|^2d\eta dx'dt\\
\ns\ds \hspace{6.1cm} \leq l^2 \int_0^{t_0} \int_{G} \theta^2 |w(t,\eta,x')|^2d\eta dx'dt.
\end{array}
\end{equation}
By virtue of $$ \n\int_0^{x_1}  w(t,\eta,x')d\eta = \int_0^{x_1} \n w(t,\eta,x')d\eta  +w(t,0,x') =  \int_0^{x_1} \n w(t,\eta,x')d\eta,  $$
 we get that
\begin{equation}\label{inv th2 ine3}
\begin{array}{ll}\ds
\int_0^{t_0} \int_G \theta^2 \Big|  \n\int_0^{x_1} w(t,\eta,x')d\eta  \Big|^2 dxdt =\int_0^{t_0} \int_G \theta^2 \Big|  \n\int_0^{x_1} w(t,\eta,x')d\eta  \Big|^2 dxdt\\
\ns\ds \hspace{6.35cm}
\leq l \int_0^l dx_1 \int_0^{t_0} \int_{G'} \int_0^l \theta^2 |\n w(t,\eta,x')|^2d\eta dx'dt\\
\ns\ds \hspace{6.35cm} \leq l^2 \int_0^{t_0} \int_{G} \theta^2 |\n w(t,\eta,x')|^2d\eta dx'dt.
\end{array}
\end{equation}
From inequality \eqref{inv th2 ine1} -- \eqref{inv th2 ine3}, we obtain that
\begin{equation}\label{inv th2 ine4}
\begin{array}{lll}\ds
\q\mathbb{E}\int_0^{t_0}\int_G \theta^2 \big(\l |\n w|^2 + s\l^2 w^2  \big)dxdt \\
\ns\ds \leq
C\mathbb{E}\int_0^{t_0}\int_{G}\theta^{2}
|\chi' u|^2 dxdt + C l^2 r_1\mathbb{E}
\int_0^{t_0}\int_G \theta^2 \Big(  |\n w |^2 +
|w|^2 \Big)dxdt.
\end{array}
\end{equation}
Thus, we know that there is a $\l_4= \max\{C
r_1,\l_1\}$ such that for all $\l\geq\l_4$,
there exists an $s_1(\l_4)>0$ so that for all $s\geq
s_1(\l_4)$, it holds that
\begin{equation}\label{inv th2 ine5}
\mathbb{E}\int_0^{t_0}\int_G \theta^2 \big(\l |\n w|^2 + s\l^2 w^2  \big)dxdt
 \leq C\mathbb{E}\int_0^{t_0}\int_{G}\theta^{2} |\chi' u|^2 dxdt.
 \end{equation}
Fix $\l = \l_4$, by the property of $\chi$(see \eqref{chi1}), we find
\begin{equation}\label{inv th2 ine6}
 \mathbb{E}\int_0^{t_0}\int_G \theta^{2} |\chi' u|^2 dxdt \leq e^{2s e^{-\l_4 t_1}} \mathbb{E}\int_Q |u|^2dxdt \leq e^{2s e^{-\l_4 t_1}} |y_{x_1}|^2_{L^2_{\cF}(0,T;L^2(G))}.
\end{equation}
This, together with inequality \eqref{inv th2
ine5}, implies that for all $s\geq s_1$, it
holds that
\begin{equation}\label{inv th2 eq7}
\begin{array}{ll}\ds
e^{2s e^{-\l_4(t_0 -\e)}} \mathbb{E}\int_0^{t_0 -\e} \int_G \big( |\n w|^2 + s w^2 \big)dxdt
  \leq  \mathbb{E}\int_0^{t_0 -\e}\int_G \theta^2 \big( |\n w|^2 + s w^2 \big)dxdt\\
 \ns\ds \hspace{7.5cm} \leq \mathbb{E}\int_0^{t_0 } \int_G \theta^2 \big( |\n w|^2 + s w^2 \big)dxdt \\
 \ns\ds \hspace{7.5cm} \leq C e^{2s e^{-\l_4 t_1}}
 |y_{x_1}|^2_{L^2_{\cF}(0,T;L^2(G))}.
\end{array}
\end{equation}
 From inequality \eqref{inv th2 eq7}, we have
\begin{equation}\label{inv th2 eq8}
|w|^2_{L^2_{\cF}(0,T;H^1(G))} \leq C e^{2s( e^{-\l_4 t_1}- e^{-\l_4(t_0 -\e)}  )} |y_{x_1}|^2_{L^2_{\cF}(0,T;L^2(G))}.
\end{equation}
Recalling that $t_0 - \e < t_1$, we know   $ e^{-\l_4 t_1}- e^{-\l_4(t_0 -\e)}<0$. Letting $s\to +\infty$, we obtain that
$$
w=0 \;\mbox{ in }\; (0,t_0 - \e)\t G, \;\; P\mbox{-a.s.}
$$
This, together with equality \eqref{inv th2 eq1},  implies that
$$
z=0 \;\mbox{ in }\; (0,t_0 - \e)\t G, \;\; P\mbox{-a.s.},
$$
which means
$$
h=0 \;\mbox{ in }\; (0,t_0 - \e)\t G', \;\;
P\mbox{-a.s.}
$$
Since $\e>0$ is arbitrary, the proof of
Theorem \ref{inv th2} is completed. \endpf

\section*{Acknowledgments.} The author  would like to thank
the anonymous referees for helpful comments.


\end{document}